\documentclass[letterpaper, 10 pt, conference]{ieeeconf}  
\usepackage{style}

\title{\LARGE \bf
Conformal Contraction for Robust Nonlinear Control with Distribution-Free Uncertainty Quantification
}
\author{Sihang Wei, Melkior Ornik, and Hiroyasu Tsukamoto
\thanks{The authors are with the Department of Aerospace Engineering, The Grainger College of Engineering, University of Illinois Urbana-Champaign, Urbana, Illinois  61801, {\tt\small sihangw2@illinois.edu, mornik@illinois.edu, hiroyasu@illinois.edu}.}%
}

\begin{document}
\maketitle
\thispagestyle{IEEEtitlestyle}
\begin{abstract}

We present a novel robust control framework for continuous-time, perturbed nonlinear dynamical systems with uncertainty that depends nonlinearly on both the state and control inputs. Unlike conventional approaches that impose structural assumptions on the uncertainty, our framework enhances contraction-based robust control with data-driven uncertainty prediction, remaining agnostic to the models of the uncertainty and predictor. We statistically quantify how reliably the contraction conditions are satisfied under dynamics with uncertainty via conformal prediction, thereby obtaining a distribution-free and finite-time probabilistic guarantee for exponential boundedness of the trajectory tracking error. We further propose the probabilistically robust control invariant (PRCI) tube for distributionally robust motion planning, within which the perturbed system trajectories are guaranteed to stay with a finite probability, without explicit knowledge of the uncertainty model. Numerical simulations validate the effectiveness of the proposed robust control framework and the performance of the PRCI tube.




\end{abstract}
\section{Introduction}
\label{sec_intro}


Robust control design for nonlinear dynamical systems with formal guarantees remains a critical challenge, especially in the presence of system uncertainty with no structural assumptions. With increased levels of autonomy and complexity,  modern control systems are increasingly deployed in safety-critical situations. However, establishing the theoretical limit of achievable guarantees with robust nonlinear control is almost essential. We approach this problem from the perspective of contraction theory and conformal prediction for a distribution-free analysis of incremental exponential boundedness under uncertainty.

Conventional control methods for uncertain systems often rely on a priori knowledge of how uncertainty and external disturbances affect the system dynamics in obtaining performance guarantees. Typically, robust control~\cite{zhou1998essentials,Khalil:1173048} assumes known bounds on worst-case disturbances or model uncertainty, stochastic control~\cite{sto_stability_book} assumes known uncertainty distributions over random disturbances or modeling errors, and adaptive control~\cite{Ref_Slotine,ioannou1996robust} assumes affine structures in how unknown parameters affects the system dynamics. Such a limitation is also present in contraction theory---a differential framework to analyze nonlinear dynamical systems by means of a contraction metric, the existence of which is a necessary and sufficient characterization of incremental exponential stability~\cite{Ref:contraction1}, it does not directly address uncertainty without additional modeling assumptions. Existing methods primarily address control design with bounded disturbances~\cite{manchester2017control,7989693}, stochastic uncertainty with a known distribution~\cite{mypaperTAC}, and structured parametric uncertainty~\cite{lopez2021universal} for formal guarantees on the exponential trajectory tracking performance. Even with data-driven approaches (see~\cite{tutorial} and the references therein), some structural assumptions on the uncertainty prediction model are still required in order to establish formal guarantees---these cannot be entirely eliminated.

Conformal prediction~\cite{conformal1, conformalbook} has recently gained attention in control theory as a powerful tool for quantifying the system uncertainty solely from past observations and without requiring explicit knowledge of the underlying distribution. Given an uncertainty predictor and a negatively oriented score function to measure the predictor's performance, it exploits the rank statistics of the score to construct a high-confidence bound on the prediction error. This concept has been successfully utilized for formal verification in data-driven control across various scenarios~\cite{confcont1}, including, but not limited to, robust motion planning~\cite{confcont2,confcont3}, adaptive motion planning~\cite{adaptive_conformal}, and stochastic optimal control of linear systems~\cite{confcont4}. However, existing work has largely focused on discrete-time settings and has not addressed continuous-time, closed-loop nonlinear control problems through contraction-based analysis. In particular, these are not compatible with analyzing the exponential nature of trajectory tracking, which is central in explicitly evaluating transient specifications of nonlinear systems. The main focus of this paper is to explore the potential of conformal prediction for distribution-free, formal verification in continuous-time robust nonlinear control, particularly through the lens of contraction theory.

\begin{figure}[t]  
    \centering
    \includegraphics[width=0.9\linewidth]{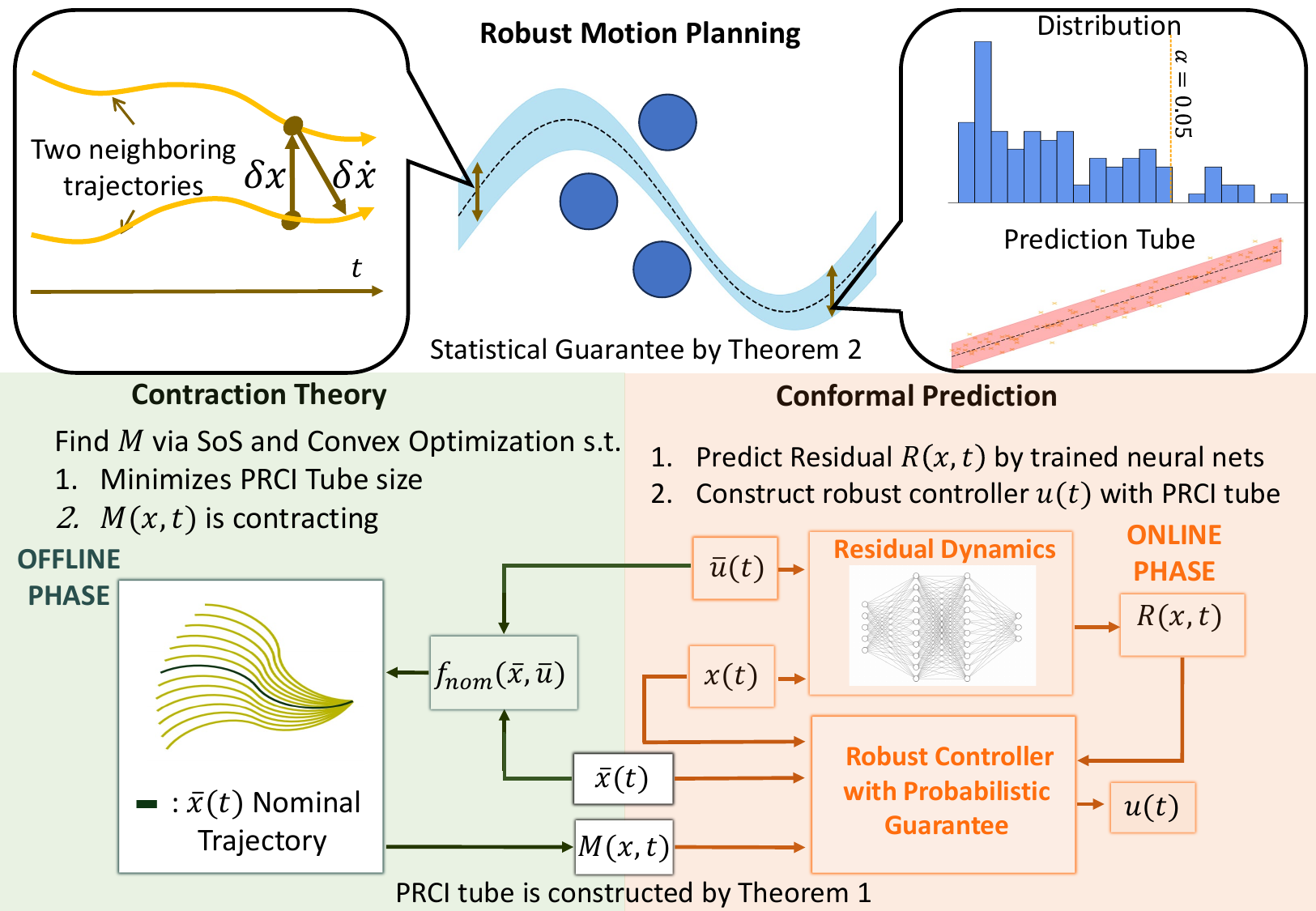}  
    \caption{From Contraction Theory to Distribution-Free Robust Nonlinear Control and Motion Planning.}
    \label{fig:intro}
\end{figure}
In this paper,we propose a novel contraction theory-based robust control framework for nonlinear dynamical systems subject to uncertainty that is nonlinear both in the state and control inputs. Our framework incorporates a data-driven uncertainty prediction model, where the discrepancy between the true (\ie, perturbed) model and the nominal model of the dynamical system can be quantified via conformal prediction as shown in Fig.~\ref{fig:intro}. This discrepancy characterizes the uncertainty of the system. 
Our main contribution is to provide a distribution-free, finite-time probabilistic guarantee for exponential boundedness of the system trajectory tracking errors, based on statistically valid satisfaction of contraction constraints under quantified uncertainty. It is worth noting that the guarantees of our proposed method remain valid regardless of the choice of the uncertainty prediction model used. We further propose the probability robust control invariant tube centered around a reference trajectory for distributionally robust motion planning, which certifies the perturbed system trajectories remain within the tube despite unstructured uncertainty.

We validate the performance of our approach through numerical simulations. Our results indicate that our proposed approach indeed achieves exponential trajectory tracking with provable probabilistic bounds on the tracking error, even under the presence of the unstructured uncertainty nonlinear both in the
state and control inputs.
\subsection*{Notation}
\label{notation}
For $A \in \mathbb{R}^{n \times n}$, we use $A \succ 0$, $A \succeq 0$, $A \prec 0$, and $A \preceq 0$ for the positive definite, positive semi-definite, negative definite, negative semi-definite matrices, respectively. For $x \in \mathbb{R}^n$, we let $\|x\|$ denote the Euclidean norm.
\section{Preliminaries and Problem Formulation}
\label{sec_prelim}
\subsection{Contraction Theory}
We present a brief review of contraction theory~\cite{Ref:contraction1,manchester2017control}, which is the fundamental tool in establishing the formal guarantees of this paper, simplifying and generalizing Lyapunov theory. Consider the following perturbed, control-affine nonlinear dynamical system
\begin{align}
    \dot{x}(t) = f(x(t)) + B(x(t)) u(x(t),t) + \zeta(x(t), u(t)),
\label{eq:sys_dyn}
\end{align}
where \(t \in \mathbb{R}_{\geq 0}\) is time, \(x:\mathbb{R}_{\geq 0} \mapsto \mathcal{X}\subseteq\mathbb{R}^n\), \(x\)  is the system state, \(u:\mathbb{R}^n\times\mathbb{R}_{\geq 0} \mapsto \mathcal{U}\subseteq\mathbb{R}^m\) is the system control policy,  ${f}:
\mathbb{R}^n\mapsto\mathbb{R}^n$ and $B:\mathbb{R}^n\mapsto\mathbb{R}^{n\times m}$ are known smooth functions, and \( \zeta : \mathbb{R}^n \times \mathbb{R}^m \mapsto \mathbb{R}^n \) is a smooth, state and control non-affine function representing the system uncertainty. Using the concepts of the finite-time stability in~\cite{1965_Weiss_FT_Stability} and the exponential boundedness in~\cite[pp. 168-174]{Khalil:1173048}, let us first define the incremental exponential boundedness under uncertainty as follows~\cite{Ref:contraction1,7989693} to describe how the distance between a perturbed and a nominal trajectory shrinks over time, up to a bounded error.
\begin{definition}[Finite-time IEB]
    Consider any pair of a perturbed state trajectory $x(t)$ of~\eqref{eq:sys_dyn} and an unperturbed reference state trajectory $\bar{x}(t)$, where $\bar{x}(t)$ is a state trajectory of the dynamics \eqref{eq:sys_dyn} giving reference control input $\bar{u}(t)$ and with \(\zeta(x,u)=0\). The perturbed trajectory \( x(t) \) evolves under a feedback control policy \( u(x, t) \) designed to keep it close to \( \bar{x}(t) \). Suppose that the control policy $u(x, t)$ can be constructed to have the following bound for any initial conditions \(x(0),\bar{x}(0) \in \mathcal{X}\):
    \begin{align}
    \begin{aligned}
        &\exists c_1,c_2,\lambda \in [0,\infty)\text{ \st{} } \\
        &d_{\scriptscriptstyle\mathrm{IEB}}(\bar{x}(t),x(t)) \leq c_1e^{-\lambda t}+c_2,~\forall t \in [0,T],
    \end{aligned}
        \label{eq:GEIS}
    \end{align}
    where $d_{\scriptscriptstyle\mathrm{IEB}}$ is the distance function that quantifies separation between the perturbed and unperturbed trajectories such as the Euclidean distance or a Riemannian distance induced by a metric at time $t$, and $T>0$ is the system time horizon. Then the system~\eqref{eq:sys_dyn} is \emph{incrementally exponentially bounded} (IEB) under uncertainty with respect to $u$, with the contraction rate \(\lambda\). If $T$ is finite, the system~\eqref{eq:sys_dyn} is finite-time IEB with respect to $u$.
    \label{def_GEIS}
\end{definition}

To analyze the IEB formally, we next define the following concepts~\cite{manchester2017control,7989693} from Riemannian geometry and introduce the contraction metrics. 
\begin{definition}[Riemannian distance]
\label{def_riemmanian}
Consider a smooth, state-dependent matrix function  \(M(x) \succ  0\). For any smooth time-varying curve \(c(\mu) : [0, 1] \mapsto \mathcal{X}\), the Riemannian length and energy are given as $L(c) = \int_0^1 \left\| c_\mu(s) \right\|_M ds$ and $E(c) = \int_0^1 \left\| c_\mu(s) \right\|_M^2 ds$, where \(c_\mu(s) = \partial c/\partial \mu\), and the norm \(\|v\|_M\) is defined as \(\|v\|_M^2 = v^\top M(x) v\). The \emph{Riemannian distance} between two points \(x, y \in \mathcal{X}\) is then defined as the follows:
\begin{align}
    \label{eq_rm_distance}
    d_{\scriptscriptstyle\mathrm{RM}}(x, y) = \inf_{c \in \Gamma(x, y)} L(c),
\end{align}
where \(\Gamma(x, y)\) is the set of all smooth curves on \(\mathcal{X}\) satisfying \(c(0) = x\) and \(c(1) = y\).
The (non-unique) geodesic that achieves the infimum is said to be the minimizing geodesic \(\gamma\). Note that $d_{\scriptscriptstyle\mathrm{RM}}(x, y)^2 = E(\gamma) = L(\gamma)^2 \leq L(c)^2 \leq E(c)$.
\end{definition}

The differential dynamics of~\eqref{eq:sys_dyn} with $\zeta(x,u)=0$ can be used to analyze the IEB of Definition~\ref{def_GEIS}. This dynamics describe how infinitesimal displacements between trajectories of \(x(t)\)  and \(\bar{x}(t)\) behave, and is given as follows:
\begin{align}
    \dot{\delta}_x =  \left(\frac{\partial f}{\partial x} + \sum_{j=1}^{m} u_j \frac{\partial b_j}{\partial x} \right) \delta_x + B(x)\delta_{u},
\label{eq:diff_ode}
\end{align}
where \(\delta_x \in T_x \mathcal{X}\) is a tangent vector to a smooth path of the states \(x\), \(\delta_{u} \in T_{u} \mathcal{U}\) is a tangent vector to a smooth path of the control inputs at \(u\), \(b_j\) is the \(j\)-th column of \(B(x,t)\), and \(u_j\)  is the \(j\)-th element of control input vector \(u\). We are now ready to introduce the following lemma of contraction theory~\cite{Ref:contraction1,manchester2017control,7989693}.
\begin{lemma}
Suppose the nominal system \eqref{eq:sys_dyn} with $\zeta(x,u)=0$ satisfies the following for all $\forall x, u, t$:
\begin{subequations}
\label{eq_allccm}
    \begin{align}
        &\underline{m} I \preceq M \preceq \overline{m} I,
        \label{eq:self_bound}\\
        &\frac{\partial b_j}{\partial x}^\top M + M \frac{\partial b_j}{\partial x} + \partial_{b_j} M = 0, \quad \forall j = 1,2,\dots,m,
        \label{eq:kill_vec}\\
        &\delta_x^\top \left(\frac{\partial f}{\partial x}^\top M + M \frac{\partial f}{\partial x} + \partial_f M + \frac{\partial M}{\partial t} + 2\lambda M \right)\delta_x < 0,
        \label{eq:cond_ccm}
    \end{align}
\end{subequations}
where \(M(x)\succ 0\), \(\lambda > 0\), \(\delta_x\) satisfies \(\delta_x^\top M B = 0\) for any \(\delta_x \neq 0\), $\partial_{p} F = \sum_{k=1}^n(\partial F/\partial x_k)p_k$ for some $p(x)\in\mathbb{R}^n$ and $F(x) \in \mathbb{R}^{n\times n}$, and the arguments are omitted for notational simplicity. Let the reference state and control trajectories be \((\bar{x}, \bar{u})\) for the nominal system of \eqref{eq:sys_dyn} with $\zeta(x,u)=0$, and let \(\gamma(\mu, t)\) be the minimizing geodesic connecting the perturbed trajectory $x(t)$ of~\eqref{eq:sys_dyn} and~$\bar{x}(t)$, such that \(\gamma(0, t) = \overline{x}(t)\) and \(\gamma(1, t) = x(t)\). If $\sup_{x,u}\|\zeta(x, u)\|_2$ is bounded, there exists a control policy $u$ that makes~\eqref{eq:sys_dyn} IEB as per Definition~\ref{def_GEIS}, and the bound is given as follows:
\begin{align}
   d_{\scriptscriptstyle\mathrm{RM}}(\bar{x}(t),x(t)) \leq d_{\scriptscriptstyle\mathrm{RM}}(\bar{x}(0),x(0))e^{-\lambda t}+\frac{\sqrt{\overline{m}}\bar{\zeta}}{\lambda},
\end{align}
where $\bar{\zeta} = \sup_{x,u}\|\zeta(x, u)\|_2$ and $d_{\scriptscriptstyle\mathrm{RM}}$ is the Riemmanian distance of Definition~\ref{def_riemmanian}.
\label{lemma_CCM}
\end{lemma}
\begin{proof}
    See~\cite{manchester2017control} and~\cite{7989693}.
\end{proof}

\subsection{Conformal Prediction}
The critical limitation of Lemma~\ref{lemma_CCM} lies in its assumption that $\sup_{x,u}\|\zeta(x, u)\|_2$ is bounded. To address this issue later in Section~\ref{sec_contribution2}, let us also present the following fundamental result of conformal prediction~\cite{conformal1, conformalbook}.
\begin{lemma}
\label{lemma:general_conformal}
Suppose that we have a prediction model $Y = \mathcal{P}(X)$ of labels $Y$ from features $X$, and consider a calibration dataset consisting of $N$ previously observed features and labels, $\mathcal{D} = \{Z^{\scriptscriptstyle(k)}\}_{k=1}^{N} =\{(X^{\scriptscriptstyle(k)}, Y^{\scriptscriptstyle(k)}) \}_{k=1}^{N}$, taking values in some space \(\mathcal{Z}\). Let \(s_{\mathcal{P}}: \mathcal{Z} \to \mathbb{R}\) be a real-valued score function, where lower values of \(s_{\mathcal{P}}(Z)\) indicate better performance of the prediction $\mathcal{P}$. 
Also, let \(s_{\mathcal{P}}^{\scriptscriptstyle(j)}\) denote the \(j\)-th smallest value among \(\{s_{\mathcal{P}}(Z^{\scriptscriptstyle(k)})\}_{k=1}^{N}\). If all the data points in $\mathcal{D}$ are exchangeable with the newly observed data point $Z^{\scriptscriptstyle(\text{NEW})} \in \mathcal{Z}$, then we have
\begin{align}
    \Pr\left[s_{\mathcal{P}}(Z^{\scriptscriptstyle(\text{NEW})}) \leq s_{\mathcal{P}}^{\scriptscriptstyle(j_{\alpha})}\right] \geq 1-\alpha,
\end{align}
where \(\alpha \in (0, 1)\) is a given miscoverage level and $j_{\alpha} = \lceil (1 - \alpha)(N + 1) \rceil$.
\end{lemma}
\begin{proof}
    See~\cite{conformal1, conformalbook}.
\end{proof}

\subsection{Problem Formulation}
\label{sec:prob}
This paper addresses the following problem.
\paragraph*{\textbf{Problem}}
Suppose that we have \textcolor{uiucblue}{(i)} a reference trajectory \((\bar{x}(t), \bar{u}(t))\) obtained by solving a finite-horizon motion planning problem for the nominal dynamical system~\eqref{eq:sys_dyn} with \(\zeta(x,u) = 0\), and \textcolor{uiucblue}{(ii)} a prediction model of the system uncertainty $\zeta(x,u)$. We aim to establish a statistically rigorous, prediction model-agnostic framework for 
\begin{enumerate}
    \item designing a control policy $u$ for the perturbed dynamical system~\eqref{eq:sys_dyn} with \(\zeta(x,u) \neq 0\), such that the closed‐loop system is finite-time IEB as in Definition~\ref{def_GEIS} with finite probability, and
    \item robustly modifying the state constraint of the motion planning problem for the reference trajectory, so that the perturbed trajectory $x(t)$ of~\eqref{eq:sys_dyn} satisfies the original state constraint also with finite probability,
\end{enumerate}
even under the presence of the system uncertainty $\zeta(x,u)$ without structural or distributional assumptions.
\section{Control Design and Uncertainty Quantification}
\label{sec_contribution1}

\subsection{Robust Control with Uncertainty Prediction}
\label{sec:u_design}
Consider a nominal contracting control policy $u_c$ for the dynamical system~\eqref{eq:sys_dyn} defined as follows:
\begin{align}
u_c(x(t),t) = \bar{u}(t) + k(x(t),\bar{x}(t))
\label{eq:feedback_ctrl}
\end{align}
where \((\bar{x}(t),\bar{u}(t))\) represents a given reference trajectory that adheres to the nominal dynamics~\eqref{eq:sys_dyn} with $\zeta(x,u) = 0$, and the contracting feedback term $k(x(t),\bar{x}(t))$ is designed with the contraction metric satisfying the conditions~\eqref{eq_allccm} with $M(x(t))$ and $\lambda$ of Lemma~\ref{lemma_CCM}, solving the following quadratic optimization problem at each time instant~\cite{primbs2000receding,manchester2017control, 7989693}:
\begin{align}
    \label{eq:qp_u}
        &k(x(t), \bar{x}(t)) = \argmin_{\kappa \in \mathbb{R}^m} \| \kappa \|^2 \\
        &\text{s.t. }-\gamma_{\mu}(1,t)^\top M(x(t)) (f(x(t))+g(x(t))(\bar{u}(t) +\kappa)) \\
        &+ \gamma_{\mu}(0,t)^\top  M(\bar{x}(t)) (f(\bar{x}(t))+g(\bar{x}(t))\bar{u}(t)) \geq\lambda E(\gamma(\cdot,t))\nonumber
\end{align}
where \(\gamma(\mu, t)\) is the minimizing geodesic connecting $x(t)$ and~$\bar{x}(t)$ with \(\gamma(0, t) = \overline{x}(t)\) and \(\gamma(1, t) = x(t)\) as in Lemma~\ref{lemma_CCM}, $\gamma_{\mu}(\mu,t) = {\partial \gamma}(\mu,t)/{\partial \mu}$, and $E(\gamma(\cdot,t))$ is the Riemannian energy as given in Definition~\ref{def_riemmanian} with $\cdot$ for explicitly indicating the integration along $\gamma$ over $\mu$. Note that the problem~\eqref{eq:qp_u} is feasible because of the contraction conditions~\eqref{eq_allccm} of Lemma~\ref{lemma_CCM}.

In order to deal with the system uncertainty $\zeta(x,u) \neq 0$, we further augment the nominal control policy~\eqref{eq:feedback_ctrl} with uncertainty prediction, resulting in our proposed control policy $u$ defined as follows:
\begin{align}
u(x(t),t) = u_c(x(t),t)-B(x(t))^\dagger \hat{\zeta}(x(t),u_{-}(t);\theta ),
\label{eq:u_closedloop}
\end{align}
where \(B(x(t))^\dagger\) is the pseudo-inverse of the actuation matrix \(B(x(t))\), \(\hat{\zeta}(x(t),u_{-}(t),t\bigr)\) is a given uncertainty predictor for $\zeta(x(t),u(t))$ with $\theta$ being the prediction hyperparameter, and $u_{-}(t) = u(x(t-\Delta t),t-\Delta t)$ for $t \geq \Delta t$ and $u_{-}(t) = 0$ otherwise for some positive constant $\Delta t > 0$. Note that $u_{-}(t)$ is introduced here to avoid the implicit dependence on $u$ in our control policy due to the control non-affine nature of the system uncertainty $\zeta(x(t),u(t))$. Practically speaking, $\Delta t$ can be regarded as the discrete sampling interval of a digital controller implemented on an underlying system that evolves in continuous time. Applying this control policy $u$ to the system~\eqref{eq:sys_dyn} gives
\begin{align}
    \dot{x} = f(x)+ B(x)u_c + \zeta(x,u) - B(x)B^{\dagger}(x)\hat{\zeta}(x,u_{-};\theta ),
\label{eq:sys_dyn_u}
\end{align}
where we have omitted some of the arguments for notational simplicity. We denote the \emph{residual error} of the closed-loop system~\eqref{eq:sys_dyn_u} as \(R(t)\), \ie{},
\begin{align}
    \label{eq_R}
    R(t) = \zeta(x(t),u(t)) - B(x(t))B(x(t))^{\dagger}\hat{\zeta}(x(t),u_{-}(t),t;\theta).
\end{align}

\subsection{Uncertainty Prediction and Quantification}
\label{sec_cp}
Suppose that we have access to datasets of \(N\) initial conditions \(x_0^{\scriptscriptstyle(k)} \in \mathcal{X}\) and nominal reference control policies \(\bar{u}^{\scriptscriptstyle(k)}: \mathbb{R}_{\geq 0} \mapsto \mathcal{U}\). We use these datasets to to integrate the nominal system~\eqref{eq:sys_dyn} with $\zeta(x,u) = 0$ for \(t\in \mathcal{T} = [0, T]\), yielding a nominal trajectory dataset \(\mathcal{D}_{\text{ref}}\) constructed as follows:
\begin{align}
\label{eq:ref_data}
\mathcal{D}_{\text{ref}} &= \left\{\left(\bar{\Phi}^{\scriptscriptstyle(k)},\bar{u}^{\scriptscriptstyle(k)}\right)\right\}_{k=1}^N
\end{align}
where $\bar{\Phi}^{\scriptscriptstyle(k)} = \{\varphi_t(x_0^{\scriptscriptstyle(k)},\bar{u}^{\scriptscriptstyle(k)})\mid t \in \mathcal{T}\}$ and $\varphi_t(x_0^{\scriptscriptstyle(k)},\bar{u}^{\scriptscriptstyle(k)})$ denotes the samples of the unperturbed state trajectories of~\eqref{eq:sys_dyn} integrated with the initial conditions $x_0^{\scriptscriptstyle(k)}$ and control policies $\bar{u}^{\scriptscriptstyle(k)}$. We partition the reference dataset $\mathcal{D}_{\text{ref}}$ into two disjoint subsets, a training subset \(\mathcal{D}_{\text{ref}}^{\text{train}}\) and a calibration subset \(\mathcal{D}_{\text{ref}}^{\text{cal}}\), such that
\begin{align}
\mathcal{D}_{\text{ref}}^{\text{train}} \cup \mathcal{D}_{\text{ref}}^{\text{cal}} = \mathcal{D}_{\text{ref}}, \quad \mathcal{D}_{\text{ref}}^{\text{train}} \cap \mathcal{D}_{\text{ref}}^{\text{cal}} = \emptyset
\end{align}
where \(|\mathcal{D}_{\text{ref}}^{\text{train}}| = N_1\), \(|\mathcal{D}_{\text{ref}}^{\text{cal}}| = N_2\), and \(N_1 + N_2 = N\).

In order to train the uncertainty predictor \( \hat{\zeta} \), we construct the following dataset \(\mathcal{D}^{\text{train}}\) by applying the nominal reference control policies \(\bar{u}^{\scriptscriptstyle(k)}\) in \(\mathcal{D}_{\text{ref}}^{\text{train}}\) to the perturbed system~\eqref{eq:sys_dyn} with \(\zeta(x,u)\neq 0\):
\begin{align}
\label{eq:train_data}
\mathcal{D}^{\text{train}} &= \left\{\left(\bar{\Psi}^{\scriptscriptstyle(k)}, \bar{u}^{\scriptscriptstyle(k)},\bar{\Zeta}^{\scriptscriptstyle(k)}\right)\right\}_{k=1}^{N_1}
\end{align}
where $\bar{\Psi}^{\scriptscriptstyle(k)} = \{\bar{\psi}_{t}^{\scriptscriptstyle(k)} \mid t \in \mathcal{T}\}$, $\bar{\Zeta}^{\scriptscriptstyle(k)} = \{\zeta(\bar{\psi}_{t}^{\scriptscriptstyle(k)},\bar{u}^{\scriptscriptstyle(k)}_t) \mid t \in \mathcal{T}\}$, $\bar{u}^{\scriptscriptstyle(k)}_t = \bar{u}^{\scriptscriptstyle(k)}(t)$, $\bar{\psi}_{t}^{\scriptscriptstyle(k)}=\psi_t(x_0^{\scriptscriptstyle(k)},\bar{u}^{\scriptscriptstyle(k)})$ denotes the training samples of the perturbed state trajectories of~\eqref{eq:sys_dyn} integrated with the initial conditions $x_0^{\scriptscriptstyle(k)}$ and control policies $\bar{u}^{\scriptscriptstyle(k)}$ of $\mathcal{D}_{\text{ref}}^{\text{train}}$, and $\zeta(\bar{\psi}_{t}^{\scriptscriptstyle(k)},\bar{u}^{\scriptscriptstyle(k)}(t))$ denotes the training samples of the uncertainties of~\eqref{eq:sys_dyn} evaluated over the sampled trajectories. We train our predictor \( \hat{\zeta}(x,u;\theta) \) to minimize the following loss function using~\eqref{eq:train_data}:
\begin{align}
\mathcal{L}(\theta) = \sum_{k=1}^{N_1} \sup_{t \in \mathcal{T}} \left\|\zeta(\bar{\psi}_{t}^{\scriptscriptstyle(k)}, \bar{u}^{\scriptscriptstyle(k)}_t) - \hat{\zeta}(\bar{\psi}_{t}^{\scriptscriptstyle(k)},\bar{u}^{\scriptscriptstyle(k)}_t;\theta)\right\|~~~~~
\label{eq:loss_fnc}
\end{align}
where $\bar{\psi}_{t}^{\scriptscriptstyle(k)}=\psi_t(x_0^{\scriptscriptstyle(k)},\bar{u}^{\scriptscriptstyle(k)})$ and $\bar{u}^{\scriptscriptstyle(k)}_t = \bar{u}^{\scriptscriptstyle(k)}(t)$ are again used for notational simplicity. It is worth mentioning that the reference trajectory \(\bar{\Phi}^{\scriptscriptstyle(k)}\) of $\mathcal{D}_{\text{ref}}^{\text{train}}$ is not directly used for training but rather for calibration as to be seen later, and
\begin{itemize}
    \item the prediction model in~\eqref{eq:loss_fnc} can be a neural network as to be demonstrated in Section~\ref{sec_simulation}, but
    \item our approach in the following works agnostically to the specific prediction method used for training.
\end{itemize}

For the calibration for conformal prediction described in Lemma~\ref{lemma:general_conformal}, we consider the control policy~\eqref{eq:u_closedloop} computed with \textcolor{uiucblue}{(i)} the sampled reference trajectories \(\mathcal{D}_{\text{ref}}^{\text{cal}} \) and \textcolor{uiucblue}{(ii)} the trained predictor $\hat{\zeta}$ that empirically minimizes $\mathcal{L}(\theta)$ of~\eqref{eq:loss_fnc} at $\theta = \theta^\ast$. We apply this control policy, denoted by $u^{\scriptscriptstyle(k)}$ for each data sample, to the perturbed system~\eqref{eq:sys_dyn} with \(\zeta(x,u)\neq 0\) to construct the following calibration dataset \(\mathcal{D}^{\text{cal}}\):
\begin{align}
\label{eq:calibration_data}
\mathcal{D}^{\text{cal}} = \left\{\left(\Psi^{\scriptscriptstyle(k)}, u^{\scriptscriptstyle(k)},\Zeta^{\scriptscriptstyle(k)}\right)\right\}_{k=1}^{N_2}
\end{align}
where $\Psi^{\scriptscriptstyle(k)} = \{\psi_{t}^{\scriptscriptstyle(k)} \mid t \in \mathcal{T}\}$, $\Zeta^{\scriptscriptstyle(k)} = \{\zeta(\psi_{t}^{\scriptscriptstyle(k)},u^{\scriptscriptstyle(k)}_t) \mid t \in \mathcal{T}\}$, $u^{\scriptscriptstyle(k)}_t = u^{\scriptscriptstyle(k)}(\psi_{t}^{\scriptscriptstyle(k)},t)$ $\psi_{t}^{\scriptscriptstyle(k)} = \psi_t(x_0^{\scriptscriptstyle(k)},u^{\scriptscriptstyle(k)})$ denotes the calibration samples of the perturbed state trajectories of~\eqref{eq:sys_dyn} integrated with $x_0^{\scriptscriptstyle(k)}$ and $u^{\scriptscriptstyle(k)}$ of $\mathcal{D}_{\text{ref}}^{\text{cal}}$ , and $\zeta(\psi_{t}^{\scriptscriptstyle(k)},u^{\scriptscriptstyle(k)}_t)$ denotes the calibration samples of the uncertainties of~\eqref{eq:sys_dyn} evaluated over the sampled trajectories. The nonconformity score for each data point in the calibration set is defined as the following maximal prediction error of the residual error of~\eqref{eq_R} over the time horizon $T$ for $k = 1,\cdots, N_2$:
\begin{align}
\label{eq:score_def}
&s_\zeta(\Psi^{\scriptscriptstyle(k)}, u^{\scriptscriptstyle(k)},\Zeta^{\scriptscriptstyle(k)}) \\
&=\sup_{t \in \mathcal{T}} \left\|\zeta(\psi_{t}^{\scriptscriptstyle(k)}, u_{t}^{\scriptscriptstyle(k)}) - B(\psi_{t}^{\scriptscriptstyle(k)}) B(\psi_{t}^{\scriptscriptstyle(k)})^{\dagger}\hat{\zeta}(\psi_{t}^{\scriptscriptstyle(k)}, u_{-}^{\scriptscriptstyle(k)};\theta^\ast)\right\|
\end{align}
where $\psi_{t}^{\scriptscriptstyle(k)}=\psi_t(x_0^{\scriptscriptstyle(k)},u^{\scriptscriptstyle(k)})$ and $u^{\scriptscriptstyle(k)}_t = u^{\scriptscriptstyle(k)}(\psi_{t}^{\scriptscriptstyle(k)},t)$ are again used for notational simplicity, $u_{-}^{\scriptscriptstyle(k)}$ is given by \(u_{-}^{\scriptscriptstyle(k)} = u_{t-\Delta t}^{\scriptscriptstyle(k)}\) for $t \geq \Delta t$ and $u_{-}(t) = 0$ otherwise, and $\theta^\ast$ is the empirically optimal hyperparameter. The system uncertainty can then be quantified as follows.
\begin{lemma}
\label{lemma_our_conformal_prediction}
Suppose that for a new test data sample of an initial state \( x_0^{\scriptscriptstyle(\text{NEW})} \) and a control policy \( u^{\scriptscriptstyle(\text{NEW})} \), the corresponding trajectory and uncertainty data point, denoted by $(\Psi^{\scriptscriptstyle(\text{NEW})}, u^{\scriptscriptstyle(\text{NEW})},\Zeta^{\scriptscriptstyle(\text{NEW})})$, is exchangeable with the ones of the calibration dataset $\mathcal{D}^{\text{cal}}$ given in~\eqref{eq:calibration_data}. Let \(s_{\zeta}^{\scriptscriptstyle(j)}\) denote the \(j\)-th smallest value among \(\{s_\zeta(\Psi^{\scriptscriptstyle(k)}, u^{\scriptscriptstyle(k)},\Zeta^{\scriptscriptstyle(k)})\}_{k=1}^{N_2}\) for $s_\zeta$ of~\eqref{eq:score_def}. Then we have the following bound:
\begin{align}
    \label{eq_our_conformal_prediction}
    \Pr\left[s_\zeta(\Psi^{\scriptscriptstyle(\text{NEW})}, u^{\scriptscriptstyle(\text{NEW})},\Zeta^{\scriptscriptstyle(\text{NEW})}) \leq s_{\zeta}^{\scriptscriptstyle(j_{\alpha})}\right] \geq 1-\alpha
\end{align}
where \(\alpha \in (0, 1)\) is a given miscoverage level and $j_{\alpha} = \lceil (1 - \alpha)(N_2 + 1) \rceil$.
\end{lemma}
\begin{proof}
This follows from Lemma~\ref{lemma:general_conformal}.
\end{proof}

Importantly, the distribution-free probabilistic bound~\eqref{eq_our_conformal_prediction} of Lemma~\ref{lemma_our_conformal_prediction} statistically rigorously quantifies
\begin{enumerate}
    \item the prediction error in $\hat{\zeta}$ for approximating the system uncertainty $\zeta$,
    \item the discretization error in $u_{-}$ for avoiding the implicit dependence on $u$ of~\eqref{eq:u_closedloop} in dealing with the control non-affine nature of $\zeta$, and
    \item the control projection error in $B(\psi_{t}^{\scriptscriptstyle(k)}) B(\psi_{t}^{\scriptscriptstyle(k)})^{\dagger}$ induced by potential underactuation.
\end{enumerate}

\section{Conformal Contraction}
\label{sec_contribution2}
\subsection{Statistical Exponential Boundedness}
\label{sec:sec_tube}
Our main result in this paper is given in the following theorem, addressing the robust control part of the problem described in Section~\ref{sec_prelim}.
\begin{theorem} 
Consider the closed-loop system~\eqref{eq:sys_dyn}, where the control input \( u \) is designed using~\eqref{eq:u_closedloop}. 
The synthesis of the controller depends on a positive definite matrix \( M(x) \) and a contraction rate \( \lambda \) as given in Lemma~\ref{lemma_CCM}, 
under the conditions~\eqref{eq:self_bound},~\eqref{eq:kill_vec}, and~\eqref{eq:cond_ccm}. 
Suppose the assumptions of Lemma~\ref{lemma_our_conformal_prediction} hold 
for the uncertainty predictor \( \hat{\zeta} \). Let \( (\Psi^{\scriptscriptstyle(\text{NEW})}, u^{\scriptscriptstyle(\text{NEW})}, \Zeta^{\scriptscriptstyle(\text{NEW})}) \) 
denote a new, unseen test sample. Assume that \( \hat{\zeta} \) was calibrated by the dataset \( \mathcal{D}^{\text{cal}} \) 
defined in~\eqref{eq:calibration_data}. We denote \( x(t) \) as the unseen state trajectory of the closed-loop system~\eqref{eq:sys_dyn} given with $(\Psi^{\scriptscriptstyle(\text{NEW})}, u^{\scriptscriptstyle(\text{NEW})},\Zeta^{\scriptscriptstyle(\text{NEW})})$ and \( \bar{x}(t) \) as the corresponding reference trajectory in~\eqref{eq:u_closedloop}. Let \( \gamma(\mu, t) \) be the minimizing geodesic connecting \( \bar{x}(t) \) and \( x(t) \), with boundary conditions \( \gamma(0, t) = \bar{x}(t)\) and \(\gamma(1, t) = x(t)\)
as defined in Definition~\ref{def_riemmanian}. Suppose that the Riemannian distance \( d_{\scriptscriptstyle\mathrm{RM}}(x(t), \bar{x}(t)) \) is locally Lipschitz. Then, it satisfies the following bound:
\begin{align}
    &\Pr \left[d_{\scriptscriptstyle\mathrm{RM}}(x(t), \bar{x}(t)) \leq \left( d_{\scriptscriptstyle\mathrm{RM}}(x(0), \bar{x}(0))- \frac{\sqrt{\overline{m}} s_{\zeta}^{\scriptscriptstyle(j_{\alpha})}}{\lambda}\right)e^{-\lambda t} \right. \nonumber \\
    &\quad \quad \quad \left.+ \frac{\sqrt{\overline{m}} s_{\zeta}^{\scriptscriptstyle(j_{\alpha})}}{\lambda}\right] \geq 1-\alpha,~\forall t\in \mathcal{T}
    \label{ineq:clf_cp}
\end{align}
where \(\alpha \in (0, 1)\) and $s^{\scriptscriptstyle(j_{\alpha})}$ are given in~\eqref{eq_our_conformal_prediction} of Lemma~\ref{lemma_our_conformal_prediction}, and $\overline{m}$ is given in~\eqref{eq:self_bound} of Lemma~\ref{lemma_CCM}. The bound~\eqref{ineq:clf_cp} implies that the system~\eqref{eq:sys_dyn} with the control policy~\eqref{eq:u_closedloop} is \emph{finite-time IEB} as per Definition~\ref{def_GEIS}, with $d_{\scriptscriptstyle\mathrm{IEB}} = d_{\scriptscriptstyle\mathrm{RM}}$ and
\begin{align}
    c_1 = \left| d_{\scriptscriptstyle\mathrm{RM}}(x(0), \bar{x}(0))- \frac{\sqrt{\overline{m}} s_{\zeta}^{\scriptscriptstyle(j_{\alpha})}}{\lambda}\right|,~c_2 = \frac{\sqrt{\overline{m}} s_{\zeta}^{\scriptscriptstyle(j_{\alpha})}}{\lambda}
\end{align}
in~\eqref{eq:GEIS}, at least with probability $1-\alpha$. 
\label{theom:d_rej}
\end{theorem}
\begin{proof}
Computing the time derivative of the Riemannian energy \(E(\gamma(\cdot,t))\) of Definition~\ref{def_riemmanian} under the feedback law \eqref{eq:u_closedloop} gives
\begin{align}
    \label{eq:e_dot}
    \begin{aligned}
    D^{\scriptscriptstyle+}{E}(\gamma(\cdot,t)) &= 2 \gamma_{\mu}(1, t)^\top M(x)\dot{x} - 2\gamma_{\mu}(0, t)^\top M(\bar{x})\dot{\bar{x}} \\
    &\leq -2\lambda E(\gamma(\cdot,t))+2\gamma_{\mu}(1, t)^\top M(x)R(t)
    \end{aligned}
\end{align}
where $D^{\scriptscriptstyle+}$ represents the upper Dini derivative, $R(t) = \zeta(x, u)-B(x)B(x)^\dagger \hat{\zeta}\ (x,u_{-};\theta)$ is the residual error given in~\eqref{eq_R}, the contraction condition of~\eqref{eq:qp_u} is used to obtain the second inequality, $\gamma_{\mu}(\mu,t) = \partial \gamma (\mu,t)/\partial \mu$, $\cdot$ is again for explicitly indicating the integration along $\gamma$ over $s$ with $\gamma$ being the minimizing geodesic, and the arguments of $E(\gamma(\cdot,t))$ are omitted for notational simplicity. Let us further factor the matrix $M$ for the Riemannian metric as \( M(x) = \Theta(x)^\top \Theta(x) \), and introduce a differential change of coordinates \(z_\mu(\mu, t) = \Theta(\gamma (\mu,t))\gamma_\mu(\mu, t)\), which transforms the tangent vector \( \gamma_s(\mu, t) \) along the geodesic into the coordinate frame induced by the metric. Since the velocity field of a geodesic is parallel to itself along the geodesic~\cite{spivak1999geometry,7989693}, \ie{}, 
\begin{align}
    \|z_\mu(\bar{\mu}, t)\|^2 = \gamma_{\mu}(\bar{\mu}, t)^\top M(\gamma (\bar{\mu},t))\gamma_{\mu}(\bar{\mu}, t) = E(\gamma(\cdot,t))
\end{align}
for any $\bar{\mu}\in[0,1]$, substituting these definitions into \eqref{eq:e_dot} yields $D^{\scriptscriptstyle+}{E}(\gamma(\cdot,t))  \leq -2\lambda E(\gamma(\cdot,t))+2\sqrt{\overline{m}}\|R(t)\|\sqrt{E(\gamma(\cdot,t))}$, resulting in
\begin{align}
       D^{\scriptscriptstyle+}d_{\scriptscriptstyle\mathrm{RM}}(x(t), \bar{x}(t))  &\leq -\lambda d_{\scriptscriptstyle\mathrm{RM}}(x(t), \bar{x}(t))+\sqrt{\overline{m}}\|R(t)\|~~~~
   \label{ineq:e_dot_cp}
\end{align}
where the relation $E(\gamma(\cdot,t)) = d_{\scriptscriptstyle\mathrm{RM}}(x(t), \bar{x}(t))^2$ is used. Applying the comparison lemma~\cite[pp. 102-103, pp. 350-353]{Khalil:1173048}, we get
\begin{align}
    &d_{\scriptscriptstyle\mathrm{RM}}(x(t), \bar{x}(t)) \\
    &\leq d_{\scriptscriptstyle\mathrm{RM}}(x(0), \bar{x}(0))e^{-\lambda t}+e^{-\lambda t}\int_{0}^{t}e^{\lambda \tau}\sqrt{\overline{m}}\|R(\tau)\|d\tau.
\end{align}
By definition of the score given in~\eqref{eq:score_def}, Lemma~\ref{lemma_our_conformal_prediction} indicates $\sup_{\tau\in\mathcal{T}}\|R(\tau)\| \leq s_{\zeta}^{\scriptscriptstyle(j_{\alpha})}$ at least with probability $1-\alpha$, which completes the proof.
\end{proof}
Let us emphasize again that the finite-time IEB result of theorem~\eqref{theom:d_rej} holds even with
\begin{enumerate}
    \item the system uncertainty $\zeta$ without structural and distributional assumptions
    \item the prediction model $\hat{\zeta}$ trained with arbitrary reasonable prediction algorithms.
\end{enumerate}
\subsection{Statistically Robust Motion Planning}
The robust motion planning part of the problem described in Section~\ref{sec_prelim} is addressed using the following concept.
\begin{definition}[PRCI tube]
Let \( \Omega : \mathcal{X} \to 2^{\mathcal{X}} \) be a set-valued mapping assigning a closed and bounded set \( \Omega(\bar{x}(t)) \subseteq \mathcal{X} \subseteq \mathbb{R}^n \) to each reference state \( \bar{x}(t) \in \mathcal{X} \) at time \( t \in \mathcal{T} \), with \( \bar{x}(t) \in \Omega(\bar{x}(t))\). 

We call the union \( \bigcup_{t \in \mathcal{T}} \Omega(\bar{x}(t)) \) a \textit{probabilistically robust control invariant} (PRCI) tube for $t\in\mathcal{T}$ if, for a given control policy \( u \), the following condition holds:
\begin{equation}
    \Pr \left( x(t) \in \Omega(\bar{x}(t)), ~\forall t \in \mathcal{T} ~\middle|~x(0) \in \Omega(\bar{x}(0)) \right) \geq 1 - \alpha,
\end{equation}
where \( \alpha \in (0,1) \). Note that this invariance must hold for all admissible realizations of the uncertainty \( \zeta(x, u) \) of~\eqref{eq:sys_dyn}. We also refer to the set \( \Omega(\bar{x}(t)) \) at a specific time as a PRCI set.
\label{def:PRCI}
\end{definition}
\begin{corollary}
\label{corollary_prci}
Consider the closed-loop system~\eqref{eq:sys_dyn} with the control policy $u$ of~\eqref{eq:u_closedloop}. If all the assumptions of Theorem~\ref{theom:d_rej} are satisfied and $d_{\scriptscriptstyle\mathrm{RM}}(x(0), \bar{x}(0)) \in[0, \sqrt{\overline{m}} s_{\zeta}^{\scriptscriptstyle(j_{\alpha})}/\lambda]$, 
a set-valued map given by
\begin{align}
    \label{eq_PRCI_tube}
    \Omega(\bar{x}(t)) = \left\{ \xi \in \mathcal{X} ~\Big|~ d_{\scriptscriptstyle\mathrm{RM}}(\xi, \bar{x}(t)) \leq \bar{d}_{\scriptscriptstyle\mathrm{RM}} = \frac{\sqrt{\overline{m}}s_{\zeta}^{\scriptscriptstyle(j_{\alpha})}}{\lambda}\right\}
\end{align}
yields a PRCI tube of Definition~\ref{def:PRCI} for $t\in\mathcal{T}$.
\end{corollary}
\begin{proof}
The statement follows from the probabilistic bound~\eqref{ineq:clf_cp} of Theorem~\ref{theom:d_rej}.
\end{proof}

\begin{figure}
    \centering
    \includegraphics[width=75mm]{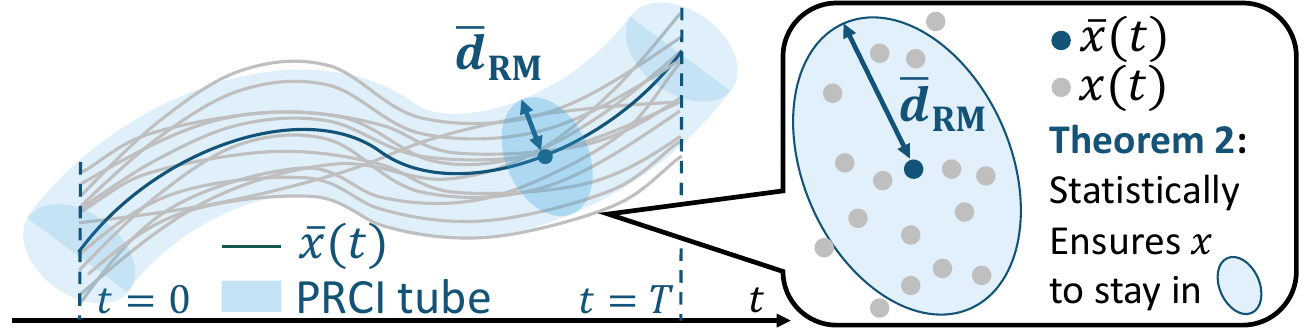}
    \caption{Illustration of the PRCI Tube in Definition~\ref{def:PRCI} and Thereom~\ref{tube_theorem}.}
    \label{fig_prci}
\end{figure}

Let us now consider the following motion planning problem for generating the reference trajectory of~\eqref{eq:u_closedloop}:
\begin{align}
\label{tube_motion_plan}
\begin{aligned}
(\bar{x},\bar{u}) &=\argmin_{s,a}~\int_{0}^{T}w_1\|a(t)\|^2+w_2P(s(t),a(t),t)dt \\
&\text{\st{} }\dot{s}(t) = f(s(t))+B(s(t))a(t),~\forall t\in\mathcal{T} \\
&\textcolor{white}{\text{\st{} }}s(t)\in\bar{\mathcal{S}}(t),~a(t)\in\bar{\mathcal{A}}(t),~\forall t\in\mathcal{T}    
\end{aligned}
\end{align}
where $w_1,w_2\geq0$, $P(s,a,t)$ is some performance-based cost function, $T$ is a given time horizon, $\mathcal{T} = [0, T]$. The sets $\bar{\mathcal{S}}(t)$ and $\bar{\mathcal{S}}(t)$ are robustly tightened admissible state and control input spaces defined as
\begin{align}
    \label{eq_tightened}
    \bar{\mathcal{S}}(t)&=\left\{\bar{\xi}\in\mathcal{X}\mid\forall \xi \in \Omega(\bar{\xi}),~\xi\in\mathcal{S}(t)\right\} \\
    \bar{\mathcal{A}}(t)&=\left\{\bar{\upsilon}\in\mathcal{U}\mid\forall \xi\in \Omega(\bar{\xi}) \text{ \& } \bar{\xi} \in \bar{\mathcal{S}}(t),~\upsilon(\xi,\bar{\xi})\in\mathcal{A}(t)\right\} \nonumber
\end{align}
where $\upsilon(\xi,\bar{\xi}) = \bar{\upsilon}+k(\xi,\bar{\xi})$ for $k$ given in~\eqref{eq:qp_u}, $\Omega$ is the PRCI mapping given in~\eqref{eq_PRCI_tube}, and $\mathcal{S}(t)$ and $\mathcal{A}(t)$ represent the original admissible state and control input spaces, respectively. The following theorem validates the effectiveness of this \textit{constraint tightening} under the presence of the unstructured uncertainty $\zeta$, as visualized in Fig.~\ref{fig_prci}.
\begin{theorem}
\label{tube_theorem}
Consider the closed-loop system~\eqref{eq:sys_dyn} with the control policy $u$ of~\eqref{eq:u_closedloop}, and suppose that all the assumptions of Corollary~\ref{corollary_prci} are satisfied. If the reference trajectories are generated by~\eqref{tube_motion_plan} with the sets~\eqref{eq_tightened} robustly tightened with the PRCI tube of Definition~\ref{def:PRCI}, then we have the following:
\begin{align}
    \Pr\left[x(t) \in \mathcal{S}(t)\text{ and } u(x(t),t) \in \mathcal{A}(t),~\forall t\in\mathcal{T}\right]\geq 1-\alpha
\end{align} 
\ie{}, the perturbed state $x(t)$ and the control input $u(x(t),t)$ both remain in the original admissible state and control input spaces, respectively, with probability at least $1-\alpha$.
\end{theorem}
\begin{proof}
Since $\bar{x}(t) \in \bar{\mathcal{S}}(t)$ and $\bar{u}(t) \in \bar{\mathcal{A}}(t)$, we have \textbf{\textcolor{uiucblue}{(i)}}~$\xi \in \mathcal{S}(t)$ for all $\xi \in \Omega(\bar{x}(t))$, and \textbf{\textcolor{uiucblue}{(ii)}}~$\bar{u}(t)+k(\xi,\bar{x}(t)) \in \mathcal{A}(t)$ for all $\xi \in \Omega(\bar{x}(t))$ by definition. Also, since the mapping $\Omega$ yields the PRCI tube as shown in Corollary~\ref{corollary_prci}, we have \textbf{\textcolor{uiucblue}{(iii)}}~$x(t) \in \Omega(\bar{x}(t))$ at least with probability $1-\alpha$ for all $t\in\mathcal{T}$. If this event occurs, viewing $\xi$ as $x(t)$ indeed implies $x(t) \in \mathcal{S}(t)$ due to \textbf{\textcolor{uiucblue}{(i)}}~and~\textbf{\textcolor{uiucblue}{(iii)}}, and $u(x(t),t) \in \mathcal{A}$ due to \textbf{\textcolor{uiucblue}{(ii)}}~and~\textbf{\textcolor{uiucblue}{(iii)}}.
\end{proof}
\begin{remark}
When we employ the constraint tightening with the calibrated PRCI tube, we restrict the test trajectories to lie within a subset of the state and input space. Doing so could change the data-generating process of Section~\ref{sec_contribution1}, potentially violating the exchangeability assumption in Lemma~\ref{lemma_our_conformal_prediction}. To address this issue, we could further split the calibration dataset into two disjoint subsets for the following two-step calibration. The first subset is used to determine the quantile bound of $\bar{d}_{\scriptscriptstyle\mathrm{RM}}$ of the PRCI tube in Corollary~\ref{corollary_prci} using Lemma~\ref{lemma:general_conformal} for its use in the constraint tightening. The second subset is then used again to perform the calibration of Lemma~\ref{lemma_our_conformal_prediction}. This minor modification allows us to maintain the coverage guarantee of Theorem~\ref{theom:d_rej} despite the potential distribution shift introduced by the constraint tightening. 
\end{remark}
\subsection{Sum of Squares Optimization for Contraction Metric}
In our numerical simulations in Section~\ref{sec_simulation}, we design the contraction metric of Lemma~\ref{lemma_CCM} with the following constraints~\cite{manchester2017control,7989693,mypaperTAC,tutorial}, shown to be equivalent to the original constraints~\eqref{eq_allccm}:
\begin{subequations}
\begin{align}
&I\preceq \bar{W}\preceq\chi I,\label{ccm_con0} \\
&B_{\bot}^{\top} \left( \partial_{b_i} \bar{W} - 2\sym{}\left({\frac{\partial b_i}{\partial x} \bar{W}}\right) \right) B_{\bot} = 0,
\label{ccm_con1} \\
&B_{\bot}^{\top} \left( -\frac{\partial \bar{W}}{\partial t}-\partial_f \bar{W} + 2\sym{}\left({\frac{\partial f}{\partial x} \bar{W}}\right) + 2 \alpha \bar{W} \right)B_{\bot} \prec 0, ~~~~\label{ccm_con2}
\end{align}
\end{subequations}
where $B_{\bot}(x)$ is a matrix whose columns span the cokernel of $B(x)$ defined as $\mathrm{coker}(B)= \{a\in\mathbb{R}^n|B^{\top}a=0\}$ satisfying $B^{\top}B_{\bot} = 0$, $b_i(x)$ is the $i$th column of $B(x)$, $\bar{W}(x)=\nu W(x)$, $W(x) = M(x)^{-1}$, $\nu=\overline{m}$, $\chi=\overline{m}/\underline{m}$, and $\partial_{p} F = \sum_{k=1}^n(\partial F/\partial x_k)p_k$ for some $p(x)\in\mathbb{R}^n$ and $F(x) \in \mathbb{R}^{n\times n}$. The convex optimization problem to minimize the squared bound of the PRCI tube in Corollary~\ref{corollary_prci} is then given as follows:
\begin{align}
    \label{eq_cvstem}
    \min_{\nu\in\mathbb{R}_{>0},\chi \in \mathbb{R},\bar{W} \succ 0, \lambda > 0}&\left(\tfrac{s_{\zeta}^{\scriptscriptstyle(j_{\alpha})}}{\lambda}\right)^2\nu\text{ \st{} \eqref{ccm_con0}, \eqref{ccm_con1}, and \eqref{ccm_con2}}
\end{align}
This problem is solved using the sum of squares programming~\cite{AYLWARD20082163} with Positivstellensatz relaxations~\cite{parrilo2003semidefinite} in Section~\ref{sec_simulation}. Note that in practice we could also consider the bound of the squared Euclidean distance given by
\begin{align}
    \|x(t)-x_d(t)\|^2 \leq \left(\tfrac{s_{\zeta}^{\scriptscriptstyle(j_{\alpha})}}{\lambda}\right)^2\tfrac{\overline{m}}{\underline{m}} = \left(\tfrac{s_{\zeta}^{\scriptscriptstyle(j_{\alpha})}}{\lambda}\right)^2\chi
\end{align}
preserving the convex structure. In particular, 
\begin{enumerate}
    \item we first solve the sum-of-squares optimization problem~\eqref{eq_cvstem} to find a feasible \( M(x) \), and
    \item we perform a line or bisection search over \(\lambda\) to identify the best trade-off between the contraction and and tube tightening.
\end{enumerate}
As future work, the ideas described in this paper could potentially be extended further to broader classes of non-Euclidean norms as in~\cite{contraction_non_euclidean,bullo_book,sontag_contraction,6632882} and stochastic nonlinear systems~\cite{Pham2009,mypaperTAC,Ref:Stochastic,han2021incremental}.

\section{Numerical Simulations}
\label{sec_simulation}
\subsection{3D Nonlinear System with Model Perturbations}
We first consider the following nominal nonlinear dynamical system~\cite{lopez2021universal} subject to parametric uncertainty and partially unmatched input disturbances, reflecting practical challenges in robust and learning-based control
\begin{align}
    \dot{x}
    =
    \begin{bmatrix} 
        x_3-\theta_1 x_1\\ 
        x_1^2 - x_2 \\ 
        \tanh(x_2) 
    \end{bmatrix}  
    + 
    \begin{bmatrix} 
        0 & 0\\ 
        1 & 0\\ 
        1 & 1
    \end{bmatrix} 
    \left(
    u - \phi(x)
    \right),
    \label{eq:ex_sys_nom}
\end{align}
where \( u = [u_1,\ u_2]^\top \),  \(\phi(x)=[\theta_2 x_3 + \theta_3 x_1^2,\ \theta_2 x_2 + \theta_3 x_1^2]^\top \) and \(\theta_1, \theta_2, \theta_3\) are constants. We assume that the true system, unknown to the planner, includes parametric uncertainties \( \delta \theta_i \) and control mismatch in the input matrix
\begin{align}
    \dot{x}
    =
    \begin{bmatrix} 
        x_3 - \tilde{\theta}_1 x_1 \\
        x_1^2 - x_2 \\
        \tanh(x_2)
    \end{bmatrix}
    +
    \begin{bmatrix} 0 & 0 \\ 0.5 & 0 \\ 1 & 0.5 \end{bmatrix} 
    (
    u - \tilde{\phi}(x)
    ),
    \label{eq:ex_sys_true}
\end{align}
where each parameter \(\theta_i\) is perturbed as \( \tilde{\theta}_i = \theta_i + \delta \theta_i \). We define \( f_{\text{nom}}(x, u) \) as the nominal dynamics in~\eqref{eq:ex_sys_nom}, and \( f_{\text{true}}(x, u) \) in~\eqref{eq:ex_sys_true} as the true dynamics. Note that the true system~\eqref{eq:ex_sys_true} can be rewritten in the additive form of~\eqref{eq:sys_dyn}, with uncertainty
\begin{align}
    \zeta(x, u) = f_{\text{true}}(x, u) - f_{\text{nom}}(x, u),
\end{align}
capturing both parametric uncertainty and control mismatch.


We consider the following setup to address the problem described in Section~\ref{sec:prob}. Specifically, our goal is to learn a data-driven uncertainty predictor \( \hat{\zeta}(x, u, t) \) to construct a PRCI tube around nominal trajectories for robust motion planning and safe trajectory tracking with a statistical guarantee under model uncertainty, as formalized in Theorems~\ref{theom:d_rej} and~\ref{tube_theorem}.

Firstly, we compute the contraction metric \(M(x)\) with the nominal values set to \(\theta_1 = 0.4, \theta_2 = 0.2, \theta_3 = 0.1\) and \(\delta \theta_1 = 0, \delta \theta_2 = 0.02, \delta \theta_3 = -0.01\). The state space is assumed to be restricted to the hypercube \([-15,\,15]^3\), and the control inputs \(u_1(t)\) and \(u_2(t)\) are limited to the interval \([-1.5,\,1.5]\). The initial states \(x(0)\) are randomly sampled uniformly within \([-15,\,15]^3\). We simulate each trajectory over a time horizon \(T=10\) with a timestep of \(\Delta t=0.01\). The uncertainty \(\zeta(x,u,t)\) is set to zero in generating the reference dataset \(\mathcal{D}_{\text{ref}}\)~\eqref{eq:ref_data}. The miscoverage level is set to \(\alpha=0.05\) for Lemma~\ref{lemma_our_conformal_prediction}. Subsequently, to train and calibrate $\hat{\zeta}(\cdot)$, we generate a dataset of 1260 trajectories for training \(\mathcal{D}^{\text{train}}\) in~\eqref{eq:train_data} and 540 trajectories for calibration \(\mathcal{D}^{\text{cal}}\) in~\eqref{eq:calibration_data} using \(u(t)\) of~\eqref{eq:feedback_ctrl}. We compare the tracking performance of a nominal controller without uncertainty compensation (\ie{},~\eqref{eq:u_closedloop} with $\hat{\zeta} = 0$) against our proposed control policy \( u(t) \) given by \eqref{eq:u_closedloop}.

\begin{figure}[htbp]  
    \centering
    \includegraphics[width=0.9\linewidth]{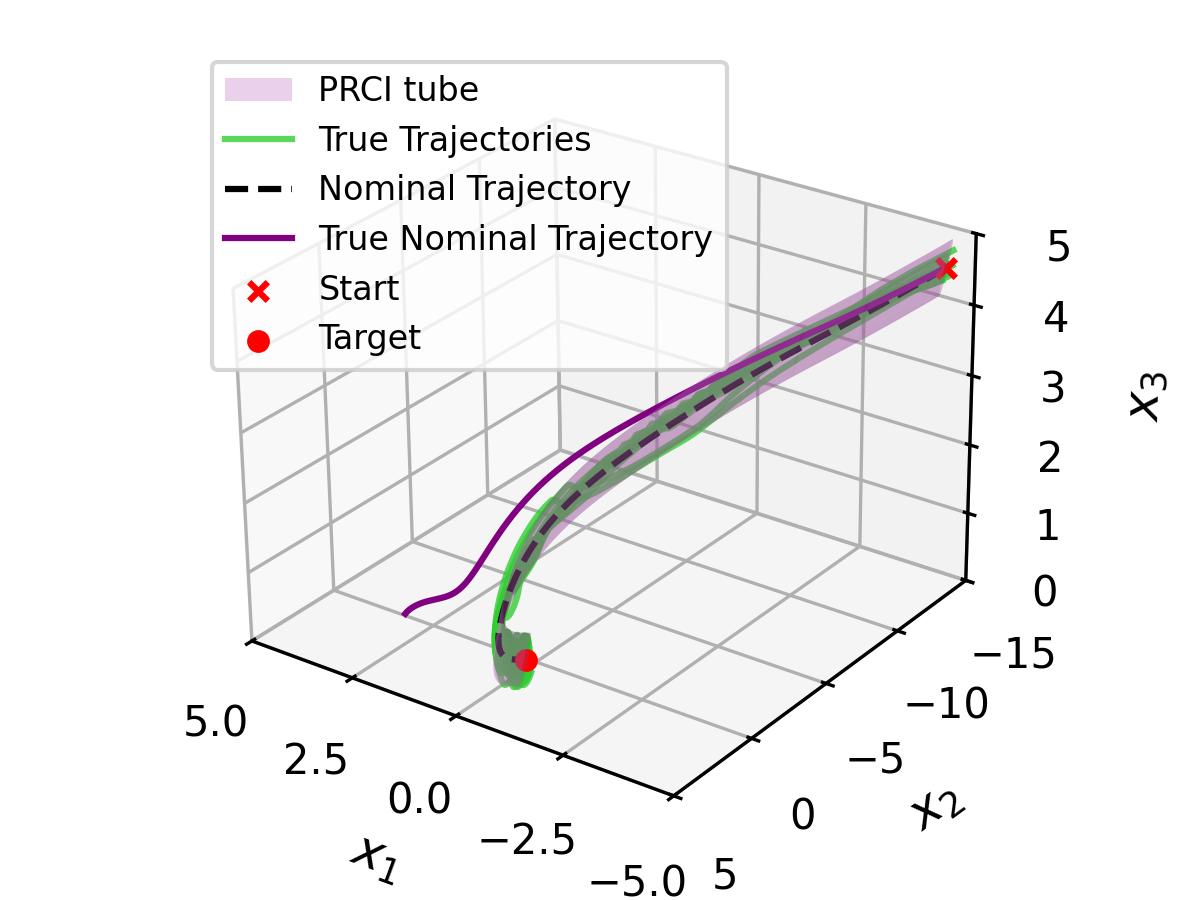}  
    \caption{Trajectories with PRCI Tube Visualization.}
    \label{fig:trajectory_tube}
\end{figure}

\begin{figure}[htbp]
    \centering
    \includegraphics[width=0.9\linewidth]{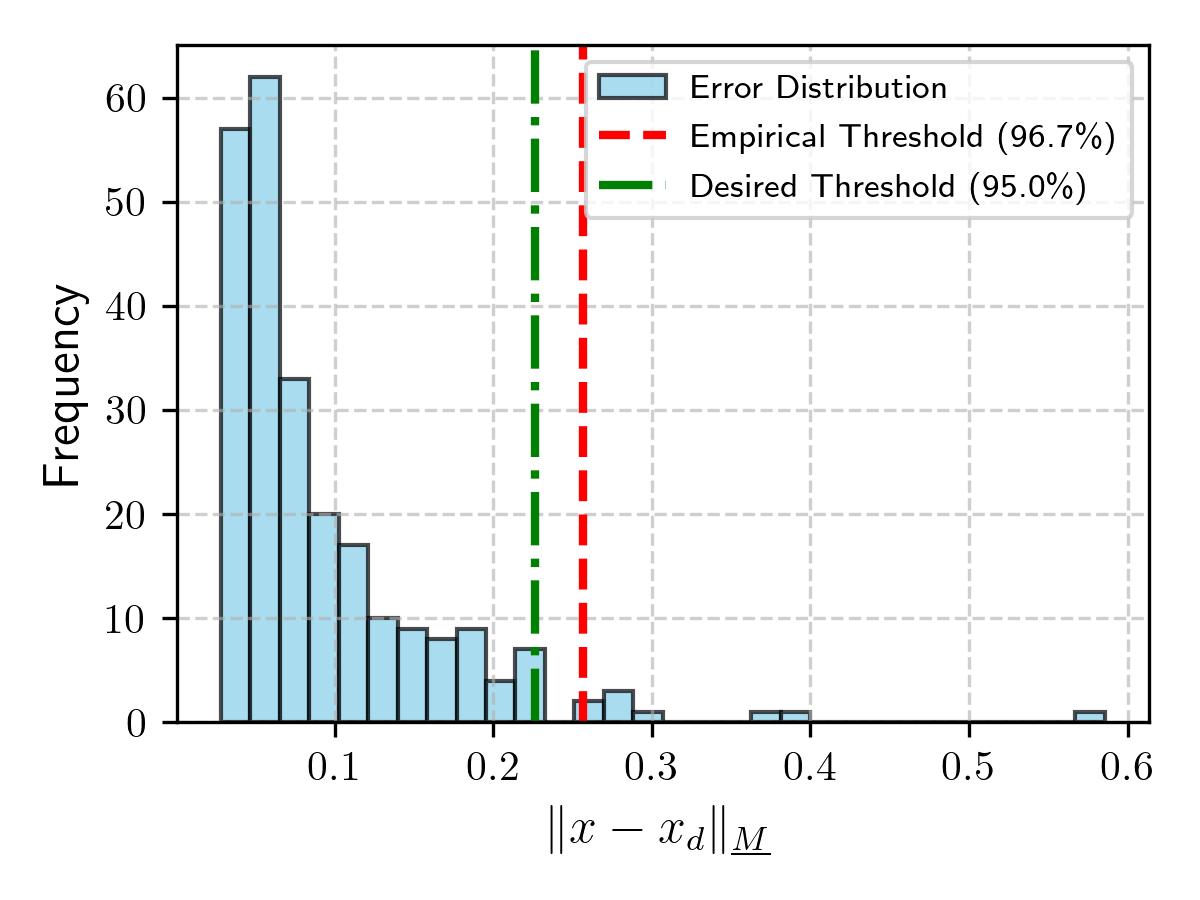}
    \caption{Empirical Distribution of Riemannian Tracking Errors.}
    \label{fig:error_exps}
\end{figure}

Figure~\ref{fig:trajectory_tube} illustrates the resulting 3D trajectory tracking performance along with the PRCI tube, where the black dashed curve represents the nominal trajectory computed under the assumption of no uncertainty, the solid green curves are the true system trajectories under the system uncertainty, controlled using our proposed control policy \( u(t) \), and the solid purple curves are the true system trajectories under the system uncertainty, controlled using the nominal control policy \(u_c(t)\). The shaded pink region represents the PRCI tube, which is a probabilistic bound within which the true trajectories are expected to remain with a high confidence level \( 1-\alpha \). To demonstrate the PRCI tube robustness guarantee of Theorems~\ref{theom:d_rej}~and~\ref{tube_theorem}, we apply our nominal control policy to 245 reference trajectories. The test trajectories remain inside the PRCI tube with empirical probability 96.7\% as shown in Fig.~\ref{fig:error_exps}, implying that our method can indeed provide a reliable and statistically calibrated bound for the worst-case deviation with the desired probability level. It can also be seen that when the system is controlled using only the nominal controller (without uncertainty rejection), most of the trajectories deviate from the PRCI tube. 
\subsection{Planar VTOL System}
As a practical example, let us also consider a planar Vertical Take-Off and Landing (VTOL) aircraft, where the state of the system is defined as \( x = [p_x, p_z, \phi, v_x, v_z, \dot{\phi}]^\top \). Here, \((p_x, p_z)\) are the horizontal and vertical positions, \(\phi\) is the roll angle, \( (v_x, v_z) \) are the lateral and vertical velocities expressed in the body-fixed frame of the vehicle, and \(\dot{\phi}\) is the angular velocity. The control input \( u = [u_1, u_2]^\top \) corresponds to the thrust produced by the two propellers. The nonlinear dynamical system is given by
\begin{align}
\dot{x} =
\begin{bsmallmatrix}
v_x \cos(\phi) - v_z \sin(\phi) \\
v_x \sin(\phi) + v_z \cos(\phi) \\
\dot{\phi} \\
v_z \dot{\phi} - g\sin(\phi) \\
-v_x \dot{\phi} - g\cos(\phi) \\
0
\end{bsmallmatrix}
+
\begin{bsmallmatrix}
0 & 0 \\
0 & 0 \\
0 & 0 \\
0 & 0 \\
1/m & 1/m \\
l/J & -l/J
\end{bsmallmatrix}
\left(
u + \zeta(x,u)
\right),
\end{align}
where \( m = 0.486\,\mathrm{kg}, J = 0.00383\,\mathrm{kg \cdot m^2}, g = 9.81\,\mathrm{m/s^2}, \) and \( l = 0.25\,\mathrm{m} \). The aircraft is tasked with navigating from an initial state set to a target set while avoiding elliptical obstacles. During flight, the system must also remain robust to the uncertainty of the model characterized by \(\zeta(x,u) = [-k_z \|v\| + k_{\dot\phi} \|u\|,\ k_{\dot\phi} \|u\|]^\top\), with \(v = [v_x, v_z]^\top\), \(k_z = 0.04\), and \(k_{\dot\phi} = 0.05\), which affects the lateral dynamics in the inertial frame.

\begin{figure}[htbp]  
    \centering
    \includegraphics[width=0.9\linewidth]{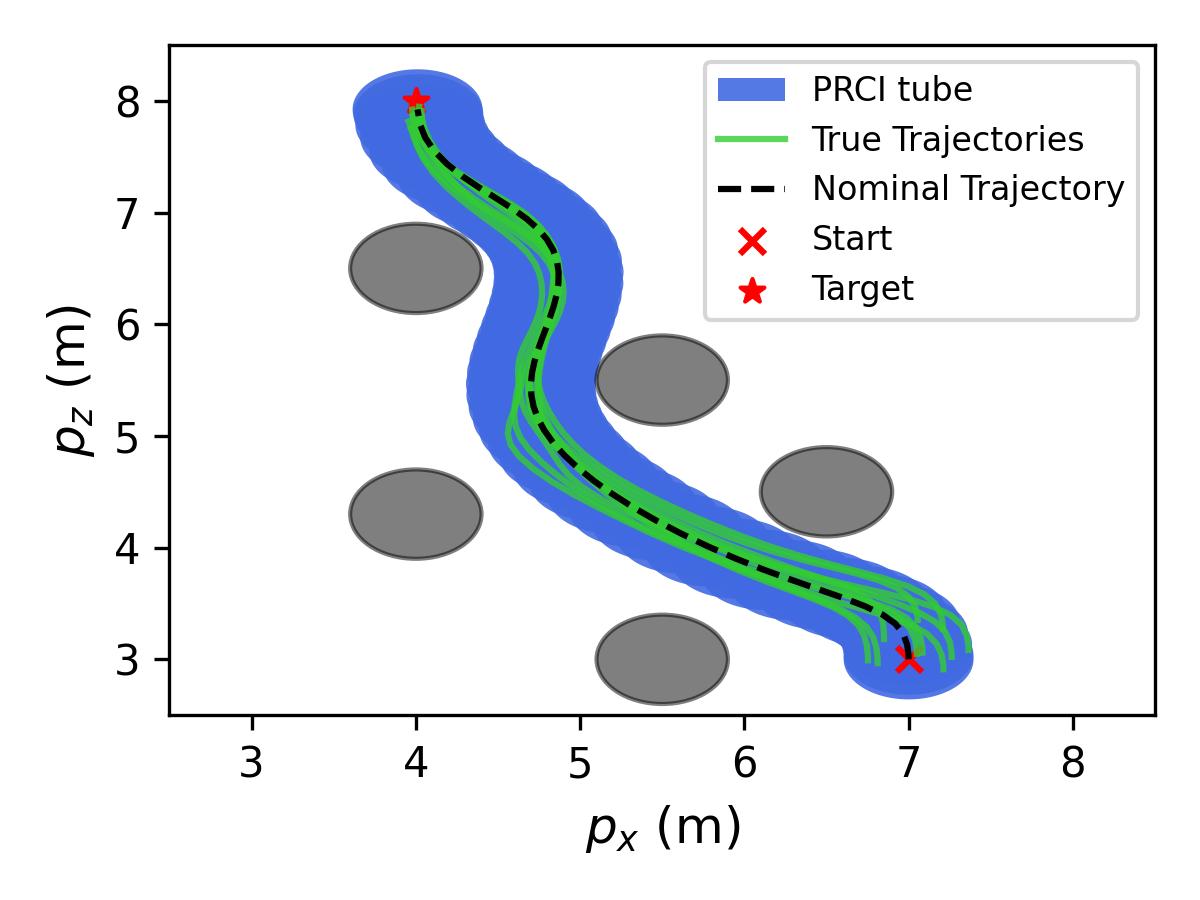}  
    \caption{Actual Trajectories with PRCI Tube in \((p_x,p_y)\) Plane Visualization.}
    \label{fig:vtol}
\end{figure}

During optimization of the contraction metric, we impose the following state bounds: \( (v_x, v_z) \in [-2,\,2] \times [-1,\,1] \) m/s, and \( (\phi, \dot{\phi}) \in [-60^\circ,\,60^\circ] \times [-60^\circ,\,60^\circ] \)/s. The miscoverage level for this PRCI tube construction is set to \(\alpha = 0.05\), which corresponds to a confidence level 95\%. We simulate each trajectory over a time horizon \(T=10\) with a timestep of \(\Delta t=0.01\). The planned nominal trajectory and multiple trajectories starting within the PRCI tube under uncertainty are illustrated in Fig.~\ref{fig:vtol}, along with the PRCI tube projected into the \(p_x, p_z\) plane by marginalizing the remaining dimensions using the Schur complement. This projection provides a 2D visualization of the high-dimensional confidence tube in the position space. The figure shows that our control policy achieves robust obstacle avoidance and safe navigation, with all executed trajectories remaining within the probabilistic bounds established by Theorems~\ref{theom:d_rej} and~\ref{tube_theorem}. In this experiment, 10 representative trajectories were tested, and no violations were observed, resulting in 100\% empirical containment. Although this number is limited due to the high computational cost of simulating long-horizon trajectories under uncertainty, the results are consistent with the target confidence level of 95\%.

\section{Conclusion}
\label{sec_conclusion}
In this work, we propose a partially data-driven framework for robust control of continuous-time nonlinear systems with uncertainty that depends nonlinearly on both the state and
control inputs. By applying an uncertainty predictor with a quantified prediction error of conformal inference, we construct a contraction theory-based control policy that achieves a finite-time, incremental exponential boundedness of Definition~\ref{def_GEIS}, without requiring strict structural assumptions on the uncertainty and its prediction models. The proposed control law is further used to construct probabilistic robust control invariant tubes for statistically verified tracking and constraint tightening in motion planning under the system uncertainty. The proposed guarantees are demonstrated using numerical simulations. Future work would include extending this approach to higher-dimensional systems, improving real-time computational efficiency, and exploring applications in online learning-based control with non-exchangeable data.
\bibliographystyle{IEEEtran}

\end{document}